\documentclass{amsart}
\usepackage{latexsym}
\usepackage{amssymb}
\usepackage{pstricks}\psset{xunit=2.2ex,yunit=2.2ex}
\newcommand{\bull}{{\bullet}}
\newcommand{\schub}{{\mathfrak S}}

\newcommand{\rr}{{\bf r}}
\theoremstyle{plain}
\newtheorem{exa}{Example}
\newtheorem{lemma}{Lemma}
\newtheorem{thm}{Theorem}
\newtheorem{prop}{Proposition}

\newcommand{\rank}{\operatorname{rank}}

\newcommand{\Schub}{{\mathfrak S}}

\newcommand{\gequ}{\geqslant}
\newcommand{\lequ}{\leqslant} 

\font\co=lcircle10
\def\petit#1{{\scriptstyle #1}}

\def\jr{\smash{\raise2pt\hbox{\co \rlap{\rlap{\char'005} \char'007}}
               \raise6pt\hbox{\rlap{\vrule height5pt}}
               \raise2pt\hbox{\rlap{\hskip4pt \vrule height0.4pt depth0pt
                width5.7pt}}
               \raise2pt\hbox{\rlap{\hskip-9.5pt \vrule height.4pt depth0pt
                width6.2pt}}
               \lower6pt\hbox{\rlap{\vrule height4.5pt}}}}
\def\je{\smash{\raise2pt\hbox{\co \rlap{\rlap{\char'005}
                \phantom{\char'007}}}\raise6pt\hbox{\rlap{\vrule height5pt}}
               \raise2pt\hbox{\rlap{\hskip-9.5pt \vrule height.4pt depth0pt
                width6.2pt}}}}
\def\er{\smash{\raise2pt\hbox{\co \rlap{\rlap{\phantom{\char'005}} \char'007}}
               \raise2pt\hbox{\rlap{\hskip4pt \vrule height0.4pt depth0pt
                width5.7pt}}
               \lower6pt\hbox{\rlap{\vrule height4.5pt}}}}
\def\+{\smash{\lower6pt\hbox{\rlap{\vrule height17pt}}
                \raise2pt\hbox{\rlap{\hskip-9pt \vrule height.4pt depth0pt
                width18.7pt}}}}
\def\hor{\smash{\raise2pt\hbox{\rlap{\hskip-9.5pt \vrule height.4pt depth0pt
                width19.2pt}}}}
\def\ver{\smash{\lower6pt\hbox{\rlap{\vrule height17pt}}}}

\def\perm#1#2{\hbox{\rlap{$\petit {#1}_{\scriptscriptstyle #2}$}}%
                \phantom{\petit 1}}

\def\mcc{\multicolumn{1}{@{}c@{}}}

\def\textcross{\ \smash{\lower4pt\hbox{\rlap{\hskip4.15pt\vrule height14pt}}
                \raise2.8pt\hbox{\rlap{\hskip-3pt \vrule height.4pt depth0pt
                width14.7pt}}}\hskip12.7pt}
\def\textelbow{\ \hskip.1pt\smash{\raise2.8pt%
                \hbox{\co \hskip 4.15pt\rlap{\rlap{\char'005} \char'007}
                \lower6.8pt\rlap{\vrule height3.5pt}
                \raise3.6pt\rlap{\vrule height3.5pt}}
                \raise2.8pt\hbox{%
                  \rlap{\hskip-7.15pt \vrule height.4pt depth0pt width3.5pt}%
                  \rlap{\hskip4.05pt \vrule height.4pt depth0pt width3.5pt}}}
                \hskip8.7pt}

\begin{document}
\title[On combinatorics of quiver component formulas]
{On Combinatorics of Quiver Component Formulas}
\author{Alexander Yong}
\date{June 2, 2003}
\address{Department of Mathematics,  
University of Michigan, 525 East University Ave.,
Ann Arbor, MI 48109-1109, USA}
\email{ayong@umich.edu}
\begin{abstract}
Buch and Fulton \cite{BF} conjectured the nonnegativity of the
{\em quiver coefficients} appearing in their formula for 
a quiver variety. Knutson, Miller and Shimozono \cite{KMS} proved this 
conjecture as an immediate consequence of their ``component formula''. We present an 
alternative proof of the component formula by substituting combinatorics for  
Gr\"{o}bner degeneration~\cite{KMS,grobGeom}. We relate the component formula
to the work of Buch, Kresch, Tamvakis and the author \cite{BKTY} where a 
``splitting'' formula for Schubert polynomials in terms of quiver coefficients was obtained. We prove
analogues of this latter result for the type $BCD$-Schubert polynomials of Billey
and Haiman~\cite{billey.haiman}.
\end{abstract}
\maketitle



\section{Introduction}

Buch and Fulton \cite{BF} established a formula for a
general kind of degeneracy locus associated to an oriented quiver of
type $A$. This formula is in terms of Schur polynomials
and certain integers, the {\em quiver coefficients}, which generalize
the classical Littlewood-Richardson coefficients. Buch and Fulton
further conjectured the nonnegativity of these quiver coefficients, and this
conjecture was recently proved by Knutson, Miller and Shimozono
\cite{KMS}. In fact, they obtained a stronger result, the ``component
formula'', whose proof was based on combinatorics, a ``ratio formula''
derived from a geometric construction due to Zelevinsky \cite{Ze} and
the method of Gr\"{o}bner degeneration, applying 
multidegree formulae for matrix Schubert varieties from \cite{grobGeom}.

In this paper, we prove a combinatorial result that
replaces the Gr\"{o}bner degeneration part of their
argument. This allows for an entirely combinatorial proof of the
component formula from the ratio formula. 
The component formula is connected to the work of Buch, Kresch, Tamvakis and the
author \cite{BKTY}, where a formula was obtained for Fulton's
universal Schubert polynomials \cite{Fulton_USP}. There, this formula was
used to obtain a ``splitting'' formula for the ordinary Schubert 
polynomials of Lascoux and Sch\"{u}tzenberger \cite{LS1} in terms of quiver coefficients. We provide 
analogues of this splitting formula for the type 
$BCD$-Schubert polynomials of Billey and
Haiman \cite{billey.haiman}, in terms of a new collection of 
positive combinatorial coefficients.

Let ${\mathfrak X}$ be a nonsingular complex variety and $E_0 \to E_1
\to \ldots \to E_n$ a sequence of vector bundles and bundle maps over
${\mathfrak X}$. A set of {\em rank conditions} for this sequence is a
collection of nonnegative integers $\rr=\{r_{ij}\}$ for $0\leq i\leq j\leq
n$. This data defines a degeneracy locus in ${\mathfrak X}$,
\[\Omega_{\rr}(E_{\bull})=\{x\in {\mathfrak X} \mid {\rm
  rank}(E_{i}(x)\to E_{j}(x))\leq r_{ij}, \forall i<j \},\]
where $r_{ii}$ is by convention the rank of the bundle $E_{i}$. We require that the
rank conditions $\rr$ {\em occur}, i.e., there exists a sequence of
vector spaces and linear maps $V_0 \to V_1 \to \cdots \to V_n$ such
that ${\rm dim}(V_i)=r_{ii}$ and ${\rm rank}(V_i\to V_j)=r_{ij}$. This
\linebreak is known to be equivalent to $r_{ij}\leq {\rm
  min}(r_{i,j-1},r_{i+1,j})$ for $i<j$ and \linebreak 
$r_{ij}-r_{i-1,j}-r_{i,j+1}+r_{i-1,j+1}\geq 0$
for $0\leq i\leq j\leq n$ where $r_{ij}=0$ if $i$ or $j$ are
not between $0$ and $n$.

The expected (and maximal) codimension of the locus
$\Omega_{\rr}(E_{\bull})$ in ${\mathfrak X}$ is 
\begin{equation}
\label{eqn:maximal_codim}
d(\rr)=\sum_{i<j} (r_{i,j-1}-r_{ij})\cdot (r_{i+1,j}-r_{ij}).
\end{equation}

Buch and Fulton \cite{BF} gave a formula for the
{\em quiver cycle}, the cohomology class of $\Omega_{\rr}(E_{\bull})$ in ${\rm H}^{*}({\mathfrak X})$,
assuming it has this codimension:
\begin{equation}
[\Omega_{\rr}(E_{\bull})]=\sum_{\underline{\mu}} c_{\underline{\mu}}(\rr) s_{\mu_{1}}(E_0 -
E_1)\cdots s_{\mu_{n}}(E_{n-1} - E_{n}).
\end{equation}
Here the sum is over sequences of partitions
$\underline{\mu}=(\mu_1,\ldots,\mu_{n})$, each $s_{\mu_{i}}$ is a super-symmetric Schur
function in the Chern roots of the bundles 
in its argument, and the \emph{quiver coefficients} $c_{\underline{\mu}}(\rr)$
are integers, conjectured to be nonnegative  by Buch and
Fulton. These coefficients generalize
the Littlewood-Richardson
coefficients, the coefficients in the expansion of a Stanley symmetric
function into Schur functions, and the
coefficients in the monomial expansions of Schubert polynomials \cite{BF, BStan, BKTY}. 

This Buch-Fulton conjecture was recently proved by Knutson, Miller and
Shimozono~\cite{KMS}. In fact, they prove the following ``component formula'':
\begin{equation}
\label{eqn:component}
[\Omega_{\rr}(E_{\bull})]=\sum_{W\in W_{\rm min}(\rr)} F_{w_{1}}(E_0 - E_1)\cdots
F_{w_n}(E_{n-1}-E_{n}),
\end{equation}
where $W_{\rm min}(\rr)$ is the set of minimal length ``lacing diagrams''
for $\rr$, and each $F_{w_{i}}$ is a double Stanley symmetric
function. 
The nonnegativity of the quiver coefficients (and a positive
combinatorial interpretation for what they count) follows immediately from (\ref{eqn:component}) by using a formula for the
expansion of a Stanley symmetric function into a positive sum of Schur functions \cite{EG,LS2}.

The set $W_{\rm min}(\rr)$ is both of combinatorial and geometric
interest. This set is derived from the strand diagrams of Abeasis and Del
Fra \cite{Abeasis_Delfra}, and generalizes the ``reduced
factorizations'' appearing in \cite{BKTY} (the latter fact is proved in
Section~5).
Moreover, this set is in 
canonical bijection with the irreducible components of 
``degenerated quiver cycles''~\cite{KMS}.

        The proof of (\ref{eqn:component}) in \cite{KMS} uses the new ``ratio formula'' for
$[\Omega_{\rr}(E_{\bull})]$, which is derived from an alternate form of a
geometric construction originally due to Zelevinsky \cite{Ze} and
developed scheme-theoretically by Lakshmibai and Magyar \cite{LM} (see
Section~\ref{subsection:combinatorics_component}) for details.
The proof proceeds by utilizing
combinatorics to derive an intermediate formula for
$[\Omega_{\rr}(E_{\bull})]$ as a multiplicity-free
sum of products of Stanley functions over {\em some} minimal length
lacing diagrams for $\rr$. Then Gr\"{o}bner geometry and Gr\"{o}bner 
degeneration \cite{KMS,grobGeom} are used to prove that {\em all} minimal length lacing
diagrams for $\rr$ actually appear.

The first goal this paper is to prove a combinatorial result that can
be substituted for the Gr\"{o}bner degeneration part of
        this proof of (\ref{eqn:component}). Combined with the rest of \cite{KMS},
this provides a combinatorial
        derivation of the component formula (\ref{eqn:component}) from
        the ratio
        formula. Our main result in this direction (Theorem 1 in Section~2)
is an explicit injection of $W_{\rm min}(\rr)$ into ${\mathcal RC}(v(\rr))$, the set of 
RC-graphs for the ``Zelevinsky permutation'' of $\rr$. This is 
proved using a characterization of Zelevinsky permutations
(Proposition~2 in Section~6).

        In Section~2, we review the definitions of some combinatorial objects associated to a
collection of rank conditions, and state our first main result. The
proof is postponed until Section~6. 
In Section~3, we 
explain the connection between this result and the proof of~(\ref{eqn:component}).
In Section~4, we provide a bijection between the labeling set in the
righthand side of  
(\ref{eqn:component}) when the rank conditions are determined by a 
permutation, and its counterpart in the formula for
Fulton's universal Schubert polynomials obtained by Buch, Kresch,
Tamvakis and the author \cite{BKTY}. 
This explains how the component formula generalizes the aforementioned
formula of~\cite{BKTY}. 

      We now describe the second goal of this paper. The
formula for the universal Schubert
polynomials obtained in \cite{BKTY} was applied there to prove a ``splitting'' formula
for the ordinary Schubert polynomials \cite{LS1} in terms of
quiver coefficients. In Section 5, we obtain
analogues of this splitting formula for the type $BCD$-Schubert polynomials
of Billey and Haiman \cite{billey.haiman}. These formulas are in terms of a 
collection of positive combinatorial coefficients that appear combinatorially 
analogous to the quiver coefficients. It would be interesting to understand
a geometric context for these formulas.
 
We thank Sergey Fomin and Ezra Miller for their questions
that initiated this work. We are extremely grateful to Ezra Miller for
introducing us to the results in \cite{KMS} and for his many helpful 
comments, including suggesting a  simplification in the proof of Theorem~1
and providing macros for drawing RC-graphs and pipe dreams.
We would also like to thank Anders Buch, Harm Derksen, Bill Fulton, 
Andrew Kresch, John Stembridge and Harry Tamvakis
for enlightening discussions. 

\section{Embedding lacing diagrams into RC-graphs}
\subsection{Ranks and laces}
Let $\rr=\{r_{ij}\}$ for $0\leq i\leq j\leq n$ be a set of rank
conditions. It is convenient to arrange them in a 
{\em rank diagram} \cite{BF}:

\[ \begin{matrix}
E_0    & \to    & E_1 & \to & E_2 & \to & \cdots & \to & E_n
\vspace{0.1cm} \\
r_{00} && r_{11} && r_{22} && \cdots && r_{nn} \\
& r_{01} && r_{12} && \cdots && r_{n-1,n} \\
&& r_{02} && \cdots && r_{n-2,n} \\
&&& \ddots \\
&&&& r_{0n}
\end{matrix} \]

We will need some notation and terminology introduced in \cite{KMS}.
The {\em lace array} $s(\rr)$ is defined by
\begin{equation}
\label{eqn:strand}
s_{ij}(\rr)=r_{ij}-r_{i-1,j}-r_{i,j+1}+r_{i-1,j+1}, 
\end{equation}
for $0\leq i\leq j\leq n$, where as before, $r_{ij}=0$ if $i$ or $j$
are not between $0$ and $n$. Note that each entry of $s_{ij}(\rr)$ is nonnegative, by our
assumptions on $\rr$.
A {\em lacing diagram} $W$ is a graph on $r_{00} + \cdots +r_{nn}$ vertices
arranged in $n$ bottom-justified columns labeled from 0 to~$n$. The
$i^{th}$ column consists of $r_{ii}$ vertices. The edges of $W$
connect consecutive columns in such a way that no two edges connecting two given columns share a vertex. A {\em lace} is a connected component of
such a graph and an $(i,j)$-{\em lace} starts in column $i$
and ends in column $j$. Also, $W$ is a {\em lacing diagram for $\rr$} if
the number of $(i,j)$-laces equals $s_{ij}(\rr)$.

\pagebreak
\begin{exa}
{\rm 
For $n=3$, the rank conditions

\setcounter{MaxMatrixCols}{7}
\[ \rr= \ \ \begin{matrix}
E_0 & \to & E_1 & \to & E_2 & \to & E_3 
\vspace{0.1cm} \\
2 && 3 && 4 && 2 \\
& 1 && 2 && 1 \\
&& 0 && 1 \\
&&& 0 \\
\end{matrix} \]
give
\[\begin{picture}(220,80)
\put(-30,30){$s(\rr)=\begin{tabular}{llll|l}
3 & 2 & 1 & 0 & $i/j$ \\ \hline
  &   &   & 1 & 0 \\ 
  &   & 0 & 1 & 1 \\ 
  & 2 & 1 & 0 & 2 \\
1 & 0 & 1 & 0 & 3 \\
\end{tabular}$} 

\thicklines
\put(120,30){$W=$}
\put(152,2){\circle*{4}}
\put(152,22){\circle*{4}}
\put(172,2){\circle*{4}}
\put(172,22){\circle*{4}}
\put(172,42){\circle*{4}}
\put(192,62){\circle*{4}}
\put(192,2){\circle*{4}}
\put(192,22){\circle*{4}}
\put(192,42){\circle*{4}}
\put(212,2){\circle*{4}}
\put(212,22){\circle*{4}}
\put(152,2){\line(1,0){20}}
\put(172,22){\line(1,-1){20}}
\put(172,42){\line(1,-1){20}}
\put(192,2){\line(1,0){20}}
\end{picture}
\]
}
\end{exa} 
\bigskip

Each lacing diagram $W$ corresponds to an ordered $n$-tuple
$(w_1,w_2,\ldots, w_n)$ of partial permutations, where $w_i$
is represented by the $r_{i-1}\times r_{i}$ $(0,1)$-matrix with an entry 1 in
position $(\alpha,\beta)$ if and only if an edge connects
the $\alpha^{th}$ vertex in column $i-1$ (counting from the bottom) to the
$\beta^{th}$ vertex in column $i$. For example, the
lacing diagram $W$ from Example~1 corresponds to:
\[\left(
\left(\begin{matrix}
1 & 0 & 0 \\
0 & 0 & 0
\end{matrix}\right), 
\left(\begin{matrix}
0 & 0 & 0 & 0 \\
1 & 0 & 0 & 0 \\
0 & 1 & 0 & 0 
\end{matrix}\right), 
\left(\begin{matrix}
1 & 0 \\
0 & 0 \\
0 & 0 \\
0 & 0  
\end{matrix}\right)\right). 
\]
Any $a\times b$ partial permutation $\rho$ has a minimal length
embedding ${\tilde \rho}$ in the symmetric group $S_{a+b}$. 
The permutation matrix
for ${\tilde \rho}$ is constructed to have $\rho$ as its northwest submatrix.
In the columns of ${\tilde \rho}$ to the right of $\rho$, place a 1 in
each of the top $a$ rows
for which $\rho$ does not already have one, making sure that the new
1's progress from northwest to southeast.  Similarly, in every row
of ${\tilde \rho}$ below ${\rho}$, place 1's going northwest to
southeast, in those columns which do not have one yet. For example, the
following are the embeddings of the above partial permutations:
\[\left(
\left(\begin{matrix}
1 & 0 & 0 & 0 & 0\\
0 & 0 & 0 & 1 & 0\\
0 & 1 & 0 & 0 & 0\\
0 & 0 & 1 & 0 & 0\\
0 & 0 & 0 & 0 & 1
\end{matrix}\right), 
\left(\begin{matrix}
0 & 0 & 0 & 0 & 1 & 0 & 0 \\
1 & 0 & 0 & 0 & 0 & 0 & 0 \\
0 & 1 & 0 & 0 & 0 & 0 & 0 \\
0 & 0 & 1 & 0 & 0 & 0 & 0 \\
0 & 0 & 0 & 1 & 0 & 0 & 0 \\
0 & 0 & 0 & 0 & 0 & 1 & 0 \\
0 & 0 & 0 & 0 & 0 & 0 & 1 
\end{matrix}\right), 
\left(\begin{matrix}
1 & 0 & 0 & 0 & 0 & 0 \\
0 & 0 & 1 & 0 & 0 & 0 \\
0 & 0 & 0 & 1 & 0 & 0 \\
0 & 0 & 0 & 0 & 1 & 0 \\
0 & 1 & 0 & 0 & 0 & 0 \\
0 & 0 & 0 & 0 & 0 & 1 
\end{matrix}\right)\right). 
\]

We define the length of a partial permutation matrix $\rho$ to be
equal to $\ell({\tilde \rho})$. Here $\ell({\tilde \rho})$ is the length
of ${\tilde \rho}$, the smallest number $\ell$ for which ${\tilde \rho}$
can be written as a product of $\ell$ simple transpositions. The
length of a lacing diagram $W=(w_1,w_2,\ldots,w_n)$ is denoted $\ell(W)$,
where $\ell(W)=\ell(w_1)+\ell(w_2)+\cdots + \ell(w_n)$. A lacing
diagram $W$ for $\rr$ is a {\em minimal length lacing diagram} if $\ell(W)=d(\rr)$. For
instance, the lacing diagram $W$ in Example~1 is of minimal length. We
denote the set of minimal length lacing diagrams for $\rr$ by $W_{\rm min}(\rr)$.
 
\subsection{The Zelevinsky permutation}
Also associated to $\rr$ is the {\em
  Zelevinsky permutation} $v(\rr)\in S_{d}$, where $d=r_{00}+r_{11}+\cdots +
  r_{nn}$ \cite{KMS}. This is defined via its {\em graph}
$G(v(\rr))$, the collection of the $d^2$ points $\{(i,w(i))\}_{1\leq i\leq
d}$ in $d\times d = [1,d]\times [1,d]$. 

Partition the $d\times d$ box into $(n+1)^{2}$ blocks $\{M_{ij}\}$
  for $0\leq
  i,j\leq n$, read as in block matrix form; so $M_{ij}$ has
dimension $r_{ii}\times r_{n-j,n-j}$ (later, we will also need the sets 
  $H_{j}=\bigcup_{i=0}^{n} M_{ij}$ and $V_{i}=\bigcup_{j=0}^{n}
  M_{ij}$ of {\em horizontal}
and {\em vertical strips} respectively).
Beginning with $M_{nn}$ and continuing right to left
  and bottom up, place $s_{n-j,i}(\rr)$ points into the block $M_{ij}$,
  as southeast as
  possible such that no two points lie in the same row or
column (in particular, points go northwest-southeast in each block). Complete the
empty rows and columns by placing points in the
super-antidiagonal blocks $M_{i,n-i-1}$, $i=0,\ldots, n-1$. In general, this
concluding step is achieved by placing points contiguously on the main
diagonal of each of the super-antidiagonal blocks. That this procedure produces
a permutation matrix is proved in \cite{KMS}.

Later we will need the fact that
\begin{equation}
\label{eqn:crucial}
\ell(v(r))=\mid\!\bigcup_{i+j\leq n-2} M_{ij}\!\mid+d(\rr).
\end{equation} 
This follows from \cite[Section~1.2]{KMS} but
can also be directly verified from (\ref{eqn:maximal_codim}) and~(\ref{eqn:strand}).

\begin{exa}\label{ex:1}{\rm 
For the rank conditions $\rr$ from Example~1, we obtain 
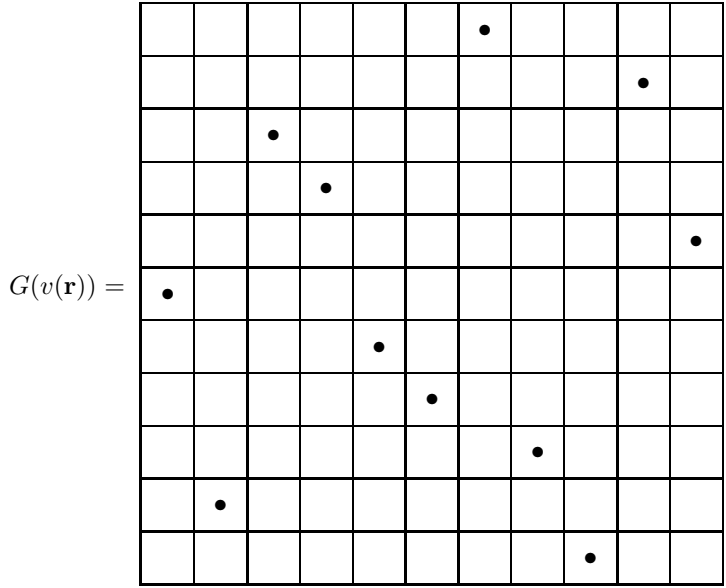
\[
\put(-50,110){$G(v(\rr))=$}
\begin{picture}(220,220)
\put(0,0){\framebox(220,220)}
\put(20,0){\line(0,1){220}}
\thicklines
\put(40,0){\line(0,1){220}}
\thinlines
\put(60,0){\line(0,1){220}}
\put(80,0){\line(0,1){220}}
\put(100,0){\line(0,1){220}}
\thicklines
\put(120,0){\line(0,1){220}}
\thinlines
\put(140,0){\line(0,1){220}}
\put(160,0){\line(0,1){220}}
\thicklines
\put(180,0){\line(0,1){220}}
\thinlines
\put(200,0){\line(0,1){220}}
\thinlines
\put(220,0){\line(0,1){220}}

\put(0,20){\line(1,0){220}}
\thicklines
\put(0,40){\line(1,0){220}}
\thinlines
\put(0,60){\line(1,0){220}}
\put(0,80){\line(1,0){220}}
\put(0,100){\line(1,0){220}}
\thicklines
\put(0,120){\line(1,0){220}}
\thinlines
\put(0,140){\line(1,0){220}}
\put(0,160){\line(1,0){220}}
\thicklines
\put(0,180){\line(1,0){220}}
\thinlines
\put(0,200){\line(1,0){220}}
\put(0,220){\line(1,0){220}}
\put(170,10){\circle*{4}}
\put(30,30){\circle*{4}}
\put(150,50){\circle*{4}}
\put(110,70){\circle*{4}}
\put(90,90){\circle*{4}}
\put(10,110){\circle*{4}}
\put(210,130){\circle*{4}}
\put(70,150){\circle*{4}}
\put(50,170){\circle*{4}}
\put(190,190){\circle*{4}}
\put(130,210){\circle*{4}}
\end{picture}
\]
Thus,
\[v(\rr)=\left(\begin{array}{cccccccccccc}
1 & 2 & 3 & 4 & 5 & 6 & 7 & 8 & 9 & 10 & 11 \\
7 & 10 & 3 & 4 & 11 & 1 & 5 & 6 & 8 & 2 & 9 
\end{array}\right).\]
}
\end{exa}

\subsection{RC-graphs}
We continue by recalling the definition of the set ${\mathcal RC}(w)$
of RC-graphs for a permutation $w\in S_{d}$ \cite{BB,FKyangBax}. 
For positive integers $a$ and $b$, consider the $a\times b$ square grid with the box in row
$i$ and column $j$ labeled $(i,j)$ as in an $a\times b$ matrix. Tile
the grid so that each box either contains a {\em cross} $\textcross$ or an 
{\em elbow joint} $\textelbow$. Thus the tiling appears as a
``network of pipes''. Such a tiled grid is a {\em pipe dream}
\cite{grobGeom}.

A {\em pipe dream for $w$} is a pipe dream where $a=b=d$, 
no crosses appear in the lower triangular part of the grid and the pipe
entering at row $i$ exits at column $w(i)$.
Finally, the set ${\mathcal RC}(w)$ of RC-graphs for a permutation $w\in
S_{d}$ is the set of pipe dreams for $w$ such that any two pipes cross
at most once. We omit drawing the ``sea of waves'' that appear
at the lower triangular part of an RC-graph.

\begin{figure}[h]
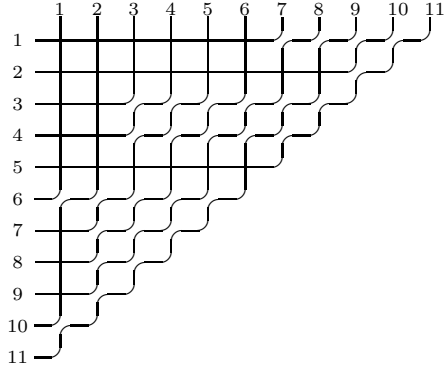

\[\begin{array}{cccccccccccc}
       &\perm1{}&\perm2{}&\perm3{}&\perm4{} & \perm5{} & \perm6{} & \perm7{} & \perm8{} 
& \perm9{} & \perm{10}{} & \perm{11}{} \\
\petit1 & \+ & \+ & \+ & \+ & \+ & \+ &  \jr  & \jr   & \jr & \jr  & \je \\
\petit2 & \+ & \+ & \+ & \+ & \+ & \+ & \+ &  \+  & \jr   & \je  &     \\
\petit3 & \+ & \+ & \jr & \jr & \jr & \jr & \jr  & \jr   &  \je &  &     \\
\petit4 & \+ & \+ & \jr & \jr & \jr & \jr &  \jr  & \je     &      &  &     \\
\petit5 & \+ & \+ & \+ & \+ & \+ & \+ &  \je  &      &      &  &     \\
\petit6 & \jr & \jr & \jr & \jr & \jr & \je &    &      &      &  &     \\
\petit7 & \+ & \jr & \jr & \jr & \je &  &    &      &      &  &     \\
\petit8 & \+ & \jr & \jr & \je &  &  &    &      &      &  &     \\
\petit9 & \+ & \jr & \je &  &  &  &    &      &      &  &     \\
\petit{10} & \jr & \je &  &  &  &  &    &      &      &  &     \\
\petit{11} & \je &  &  &  &  &  &    &      &      &  &     \\
\end{array}
\]
\caption{An RC-graph for $v(\rr)$ from Example~2}
\end{figure}

Each RC-graph is known to encode a reduced word for $w$. Let
\linebreak $u_1 u_2 \cdots u_{\ell(w)}$ be a reduced word for $w$. Then a
sequence $(\mu_{1},\mu_{2},\ldots,\mu_{\ell(w)})$ is a {\em reduced
  compatible
sequence} for $w$ if it satisfies

\begin{itemize}
\item $\mu_{1}\leq \mu_{2}\leq \ldots \leq \mu_{\ell(w)}$ 
\item $\mu_{j}\leq u_{j}$ for $1\leq j\leq \ell(w)$
\item $\mu_{j} < \mu_{j+1}$ if $u_j < u_{j+1}$ 
\end{itemize}

The following fact follows from the
definition of an RC-graph:

\begin{prop} \label{prop:rctoRC} {\rm (\cite{BB})}       
If $(\mu_{1},\ldots,\mu_{\ell(w)})$ is a reduced compatible sequence for
$w$, then
the pipe dream with crosses at $(\mu_{k},u_{k}-\mu_{k}+1)$ for $1\leq
k\leq \ell(w)$ and elbow joints elsewhere, is an RC-graph for $w$.
\end{prop}

\subsection{Main result}
Let $W=(w_1,w_2,\ldots,w_{n})$ be a lacing diagram and fix \linebreak
$\rr=\{r_{ij}\}$, $0\leq i\leq j\leq n$. A pipe dream $R$ for $w$
{\em maps to} $W$ if for all
$1\leq k\leq n$, $1\leq s\leq r_{k-1,k-1}$ and $1\leq t\leq r_{kk}$,
a pipe enters at the top of the box 
\[(r_{00}+r_{11}+\cdots+r_{k-2,k-2}+1,r_{nn}+r_{n-1,n-1}+\cdots + r_{k-1,k-1}-s+1)\]
and exits at the bottom of the box
\[(r_{00}+r_{11}+\cdots+r_{k-1,k-1},r_{nn}+r_{n-1,n-1}+\cdots + r_{kk}-t+1)\]
if and only if the $(s,t)$ entry of the partial permutation matrix
$w_k$ equals 1. Here, we set $r_{kk}=0$ if $k<0$. In other words, $R$
maps to $W$ if the above pipes correspond to the laces of $W$. 
For example, the RC-graph for $v(\rr)$ in Figure~1 maps to the lacing 
diagram $W$ from Example 1. This can be seen in the picture below: 
straightening the (partial) pipes of $W$ and right-justifying 
the result gives $W$, after reflecting across a 
northwest-southeast diagonal.

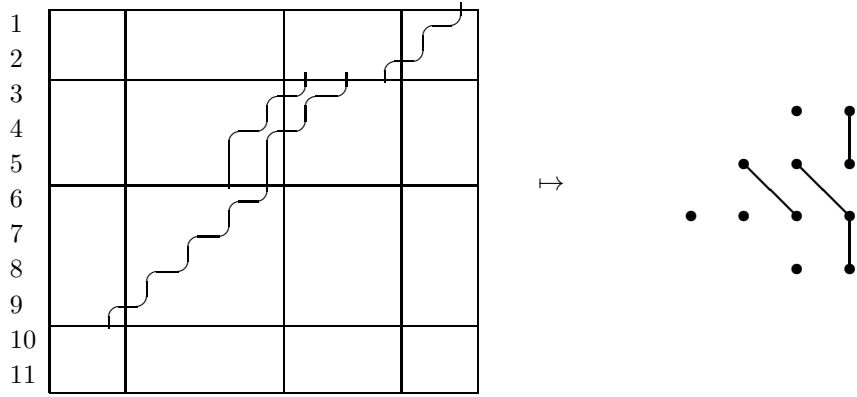
\psset{xunit=2.88ex,yunit=2.58ex}

\[
\begin{array}{@{}l|cc|cccc|ccc|cc|@{}}
\mcc{}&\ \:&\mcc{\ \,\ \,\ }&\ \:&\ \:&\ \:&\mcc{\ \,\ \,\ }&\ \:&\ \:&\mcc{\ \,\ \,\ }&\ \:&\mcc{\ \,\ \,\ }
\\[.3ex]\cline{2-12}
1  &   &   &   &   &   &   &    &    &    & \er  &\je
\\[.3ex]
2  &   &   &   &   &   &   &    &    &\er & \je  &
\\[.3ex]\cline{2-12}
3  &   &   &   &   &   & \er    &\jr&\je&   &   &
\\[.3ex]
4  &   &   &   &   &   \er         &\jr& \je&   &    &   &
\\[.3ex]
5  &   &   &   &   &\ver         & \ver&  &   &    &   &
\\[.3ex]\cline{2-12}
6  &   &   &  &    &\er&\je&   &   &    &   &
\\[.3ex]
7  &   &   &  &\er&\je&   &   &   &    &   &
\\[.3ex]
8  &   && \er &\je&   &   &   &   &    &   &
\\[.3ex]
9  &   &\er&\je&   &   &   &   &   &    &   &
\\[.3ex]\cline{2-12}
10  &   &   &   &   &   &   &   &   &    &   &
\\[.3ex]
11  &   &   &   &   &   &   &   &   &    &   &
\\\cline{2-12}
\end{array}
\quad\quad\mapsto
\begin{picture}(100,100)
\put(85,-30){\circle*{4}}
\put(105,-30){\circle*{4}}
\put(45,-10){\circle*{4}}
\put(65,-10){\circle*{4}}
\put(85,-10){\circle*{4}}
\put(105,-10){\circle*{4}}
\put(65,10){\circle*{4}}
\put(85,10){\circle*{4}}
\put(105,10){\circle*{4}}
\put(85,30){\circle*{4}}
\put(105,30){\circle*{4}}
\thicklines
\put(105,30){\line(0,-1){20}}
\put(85,10){\line(1,-1){20}}
\put(65,10){\line(1,-1){20}}
\put(105,-10){\line(0,-1){20}}
\end{picture}
\]

The following is our main result, whose proof is delayed until Section~6:
\begin{thm}
\label{thm:main_1}
Let $\rr=\{r_{ij}\}$ for $0\leq i\leq j\leq n$ be a set of rank conditions. 
There is an explicit injection from 
$W_{\rm min}(\rr) \hookrightarrow {\mathcal RC}(v(\rr))$ sending
$W\mapsto D$ such that $D$ maps to $W$.
\end{thm}
\medskip

        As explained in Section~3, combinatorics combined with the ratio
formula gives the following variation of (\ref{eqn:component}):
\[[\Omega_{\rr}(E_{\bull})]=\sum_{W\in W_{RP}(\rr)} F_{w_{1}}(E_0 - E_1)\cdots
F_{w_n}(E_{n-1}-E_{n}),\]
where $W_{RP}(\rr)$ are those $W\in W_{\rm min}(\rr)$ for which there is
a $D\in {\mathcal RC}(v(\rr))$ such that $D$ maps to $W$. Thus, 
Theorem~\ref{thm:main_1} supplies the missing ingredient for a 
combinatorial derivation of (\ref{eqn:component}) from the ratio formula.

\section{The component formula}
\subsection{Schubert polynomials}
We begin by recalling the definition of the double Schubert
polynomials of Lascoux and Sch{\"u}tzenberger \cite{L1, LS1}.
Let $X=(x_1,x_2,\ldots)$ and $Y=(y_1,y_2,\ldots)$ be two sequences of
commuting independent variables.  Given a permutation $w \in S_d$, the
double Schubert polynomial $\schub_w(X;Y)$ is defined as follows.  If
$w = w_0$ is the longest permutation in $S_d$ then we set
\[ \schub_{w_0}(X;Y) = \prod_{i+j \leq d} (x_i - y_j) \,. \]
Otherwise there is a simple transposition $s_i = (i,i+1) \in S_d$
such that \linebreak $\ell(w s_i) = \ell(w) + 1$.  We then
define
\[ \schub_w(X;Y) = \partial_i(\schub_{w s_i}(X;Y)) \]
where $\partial_i$ is the divided difference operator given by
\[ \partial_i(f) = 
  \frac{f(x_1,\dots,x_i,x_{i+1},\dots,x_d) - 
        f(x_1,\dots,x_{i+1},x_i,\dots,x_d)}{x_i - x_{i+1}}
\]
The (single) Schubert polynomial is defined by $\schub_w(X) =
\schub_w(X;0)$. By convention, if $w$ is a partial permutation, we define
$\schub_{w}=\schub_{\tilde w}$ where ${\tilde w}$ is its minimal length
embedding as a permutation.

\subsection{Symmetric functions}
Let ${\bf x}^{i}_{\bf \rr}=(x_{1}^{i},x_{2}^{i},\ldots,x_{r_{ii}}^{i})$ be
the 
Chern roots of the bundle $E_{i}$ for $0\leq i\leq n$. Then for any
partition $\lambda=(\lambda_1\geq \lambda_2\geq \ldots\geq 0)$ define
\[s_{\lambda}(E_i-E_{i+1})=s_{\lambda}({\bf x}^{i}_{\bf \rr}-{\bf
x}^{i+1}_{\bf \rr})\]
to be a super-symmetric Schur function in these roots. We will make use of the 
notation ${\bf x_{\rr}}=({\bf x}^{0}_{\rr},\ldots, {\bf x}^{n}_{\bf \rr})$
and 
${\bf {\check x_{\rr}}}=({\bf x}^{n}_{\bf \rr},\ldots, {\bf x}^{0}_{\bf \rr})$.
Similarly,
${\bf y_{\rr}}=({\bf y}^{n}_{\bf \rr},\ldots, {\bf y}^{0}_{\bf \rr})$, where
${\bf
  y_{\rr}^{i}}=(y_{1}^{i},\ldots, y_{r_{ii}}^{i})$ for $0\leq i\leq
n$. We will also need collections of infinite alphabets ${\bf x}$, ${\bf
  {\check x}}$ and ${\bf y}$, where we set $r_{ii}=\infty$ for
each $i$ in the definitions above.

For each permutation $w \in S_d$ there is a stable Schubert
polynomial or Stanley symmetric function $F_w$ in $X$
which is uniquely determined by the property that
\begin{equation}\label{Fprop}
F_w(x_1,\ldots,x_k,0,0,\ldots)=\schub_{1^m\times w}(x_1,\ldots,x_k,0,0,\ldots)
\end{equation}
for all $m\geq k$.  Here $1^m \times w \in S_{d+m}$ is the
permutation which is the identity on $\{1,\dots,m\}$ and which maps
$j$ to $w(j-m)+m$ for $j > m$ (see \cite[(7.18)]{M1}).  When $F_w$ is
written in the basis of Schur functions, one has
\begin{equation}\label{stan}
F_w=\sum_{\alpha\, : \, |\alpha|=\ell(w)} d_{w\alpha} s_{\alpha}
\end{equation}
for some nonnegative integers $d_{w\alpha}$ \cite{EG,LS2}. This also 
defines the double Stanley symmetric function $F_{w}(X-Y)$. 

\subsection{Combinatorics and the proof of (\ref{eqn:component})}
\label{subsection:combinatorics_component}

Let us now explain how our work from Section~3 leads to a
combinatorial proof of (\ref{eqn:component}). First, we
summarize the development in \cite{KMS}:

The {\em double quiver polynomial} is defined using the following {\em ratio formula}:
\[Q_{\rr}({\bf x_{\rr}};{\bf y_{\rr}})=\frac{\schub_{v(\rr)}({\bf x_{\rr}};{\bf y_{\rr}})}{\schub_{v(Hom)}({\bf x_{\rr}};{\bf y_{\rr}})},\]
where 
\[\schub_{v(Hom)}({\bf x_{\rr}};{\bf
  y_{\rr}})=\prod_{\tiny \begin{array}{c}
i+j\leq n-2\\
\alpha\leq r_{ii}, \beta\leq r_{n-j,n-j}\end{array}}
 (x_{\alpha}^{i}-y_{\beta}^{n-j})\]
It is an easy consequence of known facts about double Schubert
polynomials (see, e.g., \cite{FKyangBax}) and the definition of $v(\rr)$ that $\schub_{v(Hom)}$
divides $\schub_{v(\rr)}$.

For an integer $m\geq 0$, let $m+\rr$ be the set of rank conditions
 $\{m+r_{ij}\}$, for \linebreak $0\leq i\leq j\leq n$. It is shown that the limit
\begin{equation}
\label{eqn:stable_quiver}
F_{\rr}({\bf x}-{\bf y}):=\lim_{m\to \infty} Q_{m+\rr}({\bf x}-{\bf y})
\end{equation}
exists \cite[Proposition~6.3]{KMS}. That is, the coefficient of any fixed monomial eventually
becomes constant.

Recall $W_{RP}(\rr)$ is the set of those $W\in W_{\rm min}(\rr)$ for which there is
a $D\in {\mathcal RC}(v(\rr))$ such that $D$ maps to $W$. It is proved
combinatorially that
\begin{equation}
F_{\rr}({\bf x_{\rr}}-{\bf y_{\rr}})=\sum_{W\in W_{RP}(\rr)} F_{w_1}
({\bf x}^{0}_{\bf \rr}-{\bf
  y}^{1}_{\bf \rr})\cdots F_{w_n}({\bf x}^{n-1}_{\bf \rr}-{\bf
  y}^{n}_{\bf \rr}),
\end{equation}
where $W=(w_1,\ldots,w_n)$.

There are two facts coming from geometry that are needed.
The first is:
\begin{equation}
\label{eqn:geom1}
[\Omega_{\rr}]=Q_{\rr}({\bf x}-{\bf {\check x}})
\end{equation}
which is derived from an alternate form of a geometric construction originally due
to Zelevinsky \cite{Ze}, and developed scheme-theoretically by Lakshmibai and Magyar \cite{LM}
(and also reproved in \cite{KMS}). The second is: 
\begin{equation}
\label{eqn:geom2}
c_{\underline{\mu}}(m+\rr)=c_{\underline{\mu}}(\rr)
\end{equation}
for all $\mu$ and $m\geq 0$, which is a consequence of the main
theorem of \cite{BF}.

\medskip
\noindent
By (\ref{eqn:geom1}) and the main theorem of \cite{BF}, one has 
\[Q_{\rr}({\bf x}-{\bf {\check x}})=\sum_{\underline \mu}
c_{\underline{\mu}}(\rr)s_{\mu_1}({\bf x}^{0}_{\bf \rr}-{\bf x}^{1}_{\bf \rr})\cdots
s_{\mu_{n}}({\bf x}^{n-1}_{\bf \rr}-{\bf x}^{n}_{\bf \rr})\]

Since this holds for any ranks $\rr$, it holds for $m+\rr$ when $m$ is
large. By (\ref{eqn:stable_quiver}) and (\ref{eqn:geom2}),
\[F_{\rr}({\bf x}-{\bf {\check x}})=\sum_{\underline{\mu}}
c_{\underline{\mu}}(\rr)s_{\mu_1}({\bf x}^{0}-{\bf x}^{1})\cdots
s_{\mu_{n}}({\bf x}^{n-1}-{\bf x}^{n}).\]
Then (\ref{eqn:component}) follows after specializing ${\bf x}^{i}$ to
${\bf x_{\rr}}^{i}$ for each $i$, i.e., by setting all ``tail'' variables
$x_{j}^{i}$ for $j\geq r_{ii}+1$ to zero.

At this point, this argument gives a formula for $[\Omega_{\rr}(E_{\bull})]$ as a
multiplicity-free sum of products of Stanley functions over {\em some}
minimal length lacing diagrams for $\rr$. It remains to
show that actually {\em all} appear. The proof of this fact in
\cite{KMS} was obtained from the geometric method of Gr\"{o}bner
degeneration, by subsequently applying multidegree formulae for matrix Schubert varieties from \cite{grobGeom}. However, this is also immediate from 
Theorem~\ref{thm:main_1}. This completes a combinatorial
derivation of (\ref{eqn:component}) from the ratio formula (although we emphasize that the proof
of the latter very much depends on geometry). Note that in this proof, facts coming
from geometry are only required in order to connect the combinatorics
of the polynomials above to quiver cycles. 

In \cite{BF}, an explicit positive combinatorial formula was
conjectured for $c_{\underline{\mu}}(\rr)$. 
This is proved in \cite{KMS} using combinatorics, together with
the ratio formula and the component formula. 
Thus, Theorem~1 also allows for a combinatorial proof of that 
conjecture, starting from the ratio formula.

\section{Relations to Fulton's universal Schubert polynomials}

In this section, we report on the details of a bijection which shows how
the component formula (\ref{eqn:component}) generalizes a formula for Fulton's universal
Schubert polynomials given in \cite{BKTY}. This bijection was
also found independently in \cite{KMS}, where a proof was sketched. We
provide another proof below.

Let ${\mathfrak X}$ be a nonsingular complex variety and let
\begin{equation}
\label{eqn:univseq}
G_1 \to \cdots \to G_{n-1} \to G_n \to H_n \to
H_{n-1} \to \dots \to H_1
\end{equation}
be a sequence of vector bundles and morphisms over ${\mathfrak X}$,
such that $G_i$ and $H_i$ have rank $i$ for each $i$.
For every permutation $w$ in the symmetric group $S_{n+1}$
there is a degeneracy locus
\[
\Omega_w(G_\bull \to H_\bull)=
  \{ x \in {\mathfrak X} \ | \ \rank(G_q(x) \to H_p(x)) \lequ r_w(p,q)
  \text{ for all $1\lequ p,q\lequ n$} \},
\]
where $r_{w}(p,q)$ is the number of $i\lequ p$ such that $w(i)\lequ
q$.  The universal double Schubert polynomial $\Schub_w(c;d)$ of
Fulton \cite{Fulton_USP} gives a formula for this locus; this is a polynomial in the
Chern classes $c_i(j)=c_{i}(H_{j})$ and $d_i(j)=c_{i}(G_j)$ for $1~\lequ~i~\lequ~j\lequ~n$.
These polynomials are known to specialize to the single and double
Schubert polynomials and the quantum Schubert polynomials
\cite{FGP,C-F}.

The loci associated with universal Schubert polynomials are special   
cases of these quiver varieties.  Given $w \in S_{n+1}$ we define rank
conditions $\rr^{(n)}(w) = \{ r^{(n)}_{ij} \}$ for \linebreak $1 \lequ i \lequ j \lequ
2n$ by
\begin{equation}
\label{eqn:perm_rank_conditions}
 r^{(n)}_{ij} = \begin{cases}
   r_w(2n+1-j, i) & \text{if $i\lequ n < j$} \\
   i & \text{if $j \lequ n$} \\
   2n+1-j & \text{if $i \gequ n+1$.}
\end{cases} 
\end{equation}
The expected (and maximal) codimension of this locus is
$\ell(w)$.

Thus the quiver polynomial specializes to give a formula for the
universal
Schubert polynomial. We say that a product $u_1 \cdots u_{2n-1}$ is a
{\em reduced factorization} of $w$ if $u_1 \cdots u_{2n-1}=w$ and 
$\ell(u_1)+\cdots + \ell(u_{2n-1})=\ell(w)$.
The following was proved:\footnote{See also \cite{BKTY1} for a K-theoretic
generalization.}

\begin{thm}
\label{thm:BKTY_main}
{\rm (\cite{BKTY})}
For $w\in S_{n+1}$,
\[[\Omega_{\rr^{(n)}(w)}]=\sum_{u_1 u_2 \cdots u_{2n-1} = w}
F_{u_1}(G_1-G_2)\cdots F_{u_{2n-1}}(H_2 - H_1)\]
where the sum is over all reduced factorizations $w=u_1 \cdots
u_{2n-1}$ such that \linebreak $u_i \in S_{{\rm min}(i,2n-i)+1}$ for each~$i$.
\end{thm}

There does not appear to be any {\em a priori} reason, such as by
linear independence or geometry,
that proves that this expansion coincides
with (\ref{eqn:component}) under the conditions (\ref{eqn:univseq})
and (\ref{eqn:perm_rank_conditions}). However, this follows from:

\begin{prop}
\label{prop:bij}
The map $\Gamma$ that sends $W=(w_1,\ldots,
w_{2n-1})\in W_{\rm min}(\rr_{w}^{(n)})$ to \linebreak ${\tilde w_{2n-1}}^{-1}{\tilde
  w_{2n-2}}^{-1}\cdots {\tilde w_{1}}^{-1}$ is a bijection 
between minimal
length lacing diagrams of $\rr_{w}^{(n)}$ and reduced
factorizations of $w=u_1 \cdots u_{2n-1}$ such that $u_i\in S_{{\rm
    min}(i,2n-i)+1}$ for each~$i$.
\end{prop} 

\begin{exa}
\label{exa:bij}
{\rm 
Let $n=2$ and $w=s_2s_1=\left(
\begin{array}{ccccc}
1 & 2  & 3 \\
3 & 1  & 2 
\end{array}
\right)\in S_3.$ 
This corresponds to the following rank conditions:
\[ \rr^{(2)}(w)= \ \ \begin{matrix}
1 && 2 && 2 && 1 \\
& 1 && 1 && 1 \\
&& 1 && 0 \\
&&& 0 \\
\end{matrix} \]
The unique lacing diagram associated to $\rr^{(2)}(w)$ is
drawn below with bold lines and solid vertices. By drawing 
``phantom'' laces and vertices, $w$ is encoded by reading the 
paths from right-to-left. 

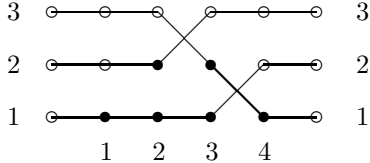
\[\begin{picture}(160,60)
\put(50,1){$1$}
\put(70,1){$2$}
\put(90,1){$3$}
\put(110,1){$4$}

\put(15,14){$1$}
\put(15,34){$2$}
\put(15,54){$3$}
\put(32,17){\circle{4}}
\put(32,37){\circle{4}}
\put(32,57){\circle{4}}
\put(32,17){\line(1,0){20}}
\put(32,37){\line(1,0){40}}
\put(32,57){\line(1,0){40}}

\thicklines
\put(52,17){\line(1,0){40}}
\thinlines
\put(72,37){\line(1,1){20}}
\thinlines
\put(72,57){\line(1,-1){20}}
\put(52,17){\circle*{4}}
\put(52,37){\circle{4}}
\put(52,57){\circle{4}}
\put(72,17){\circle*{4}}
\put(72,37){\circle*{4}}
\put(72,57){\circle{4}}

\thicklines
\put(92,37){\line(1,-1){20}}
\thinlines
\put(92,57){\line(1,0){40}}
\put(92,17){\line(1,1){20}}
\put(112,17){\line(1,0){20}}
\put(112,37){\line(1,0){20}}
\put(92,17){\circle*{4}}
\put(92,37){\circle*{4}}
\put(92,57){\circle{4}}
\put(112,17){\circle*{4}}
\put(112,37){\circle{4}}
\put(112,57){\circle{4}}

\put(132,17){\circle{4}}
\put(132,37){\circle{4}}
\put(132,57){\circle{4}}
\put(147,14){$1$}
\put(147,34){$2$}
\put(147,54){$3$}
\end{picture}
\]
}
\end{exa}

\pagebreak
\noindent
{\em Proof of Proposition~\ref{prop:bij}.}
The following lemma is an easy consequence of the definition of
$r_{w}(p,q)$:

\begin{lemma}
\label{lemma:rank}
Let $w\in S_{n+1}$, then
$r_{w}(p,q)-r_{w}(p-1,q)-r_{w}(p,q-1)+r_{w}(p-1,q-1)$ is equal to 1 if
$w(p)=q$ and is equal to 0 otherwise. Here we set $r_{w}(p,q)=0$ if
$p<0$ or $q<0$.
\end{lemma}
Lemma~\ref{lemma:rank} combined with (\ref{eqn:strand}) and
(\ref{eqn:perm_rank_conditions}) implies that
$s_{ij}(\rr_{w}^{(n)})$ for $1\leq i\leq j\leq 2n$ is 1 if $(i,j)$ falls into one of the
following three cases:
\begin{itemize}
\item[(i)] $(w(\alpha),2n-\alpha+1)$ and $1\leq w(\alpha)\leq n,
  1\leq \alpha\leq n$;
\item[(ii)] $(w(n+1),n)$ and $w(n+1)\neq n+1$;
\item[(iii)] $(n+1,2n-w^{-1}(n+1)+1)$ and $w^{-1}(n+1)\neq n+1$;
\end{itemize}
and is equal to $0$ otherwise.

First, we check that $\Gamma$ is well-defined. If 
$W=(w_1,w_2,\ldots,w_{2n-1})\in W_{\rm min}(\rr_{w}^{(n)})$ then it
is immediate from (\ref{eqn:perm_rank_conditions}) that ${\tilde w_{2n-i}}^{-1}\in S_{{\rm
    min}(i,2n-i)+1}$ for $1\leq i\leq 2n-1$. Also the conditions (i), (ii) and (iii) are
exactly saying that ${\tilde w_{2n-1}}^{-1}{\tilde w_{n-1}}^{-1}\cdots {\tilde
  w_{1}}^{-1}=w$ (e.g., by generalizing the picture in Example~\ref{exa:bij}). Further, since
\[\ell({\tilde w_{2n-1}}^{-1})+\cdots +\ell({\tilde w_{1}}^{-1})=d(\rr_{w}^{(n)})=l(w),\]
this factorization of $w$ is reduced. 

It is clear that $\Gamma$ is injective. 
To check surjectivity, let $u_1 u_2 \cdots u_{2n-1}$ be a reduced factorization of
$w$ such that $u_{i}\in S_{{\rm min}(i,2n-i)+1}$. Then let 
$W=(u_{2n-1},\ldots,u_{1})$ be the lacing diagram obtained by
interpreting each $u_{2n-i}$ as the partial permutation represented by a
${\rm min}(i,2n-i)\times ({\rm min}(i,2n-i)+1)$ matrix, for $i<n$
and a \linebreak $({\rm min}(i,2n-i)+1)\times {\rm min}(i,2n-i)$ matrix for $i>n$, and an $n\times n$ matrix for $i=n$ (in the last case, we ignore $n+1$ in
the domain and range of $u_{n}$).
This combined with $u_1\cdots u_{2n-1}=w$ shows there is a 
unique $(i,j)$-lace when one of the conditions (i),(ii) or
(iii) hold, and no other laces. Thus our calculation of $s(\rr_{w}^{(n)})$ shows 
$W$ is a lacing diagram
for $\rr_{w}^{(n)}$. This lacing diagram is of minimal length
since  $u_1 u_2 \cdots u_{2n-1}=w$ is a reduced factorization and $\ell(w)=d(\rr_{w}^{(n)})$. Finally,  \linebreak
$\Gamma$ maps $W$ to $u_1 u_2 \cdots u_{2n-1}$, as desired.\qed

\section{Splitting Schubert polynomials for classical Lie types}

        In this section, we present ``splitting'' formulas for Schubert
polynomials in each of the classical Lie types, i.e., formulas for 
polynomial representatives of Schubert classes in the cohomology
ring of generalized flag varieties \cite{BGG,borel}. In \cite{BKTY}, a splitting
formula for the Schubert polynomials of \cite{LS1} was deduced from
Theorem~\ref{thm:BKTY_main}. Our analogues use the Schubert
polynomials of types $B_n,C_n$ and $D_n$ defined by Billey and Haiman \cite{billey.haiman}.

        For a permutation
$w\in S_{n}$ and a sequence of nonnegative integers $\{a_{j}\}$ with
$1\leq a_1<a_2<\ldots<a_k<n$, we say that $w$ is {\em compatible} with 
$\{a_{j}\}$ if whenever $\ell(ws_i)<\ell(w)$ for a simple transposition
$s_i$, then $i\in\{a_{j}\}$. Also, let ${\rm col}(T)$ denote the column word
of a semi-standard Young tableau $T$, the word obtained by reading the 
entries of the columns of the tableau from bottom to top and left to right.
The following is the splitting formula for the $A_{n-1}$ Schubert
polynomials of \cite{LS1}:
\begin{thm}
\label{thm:Schub_split}
{\rm (\cite{BKTY})} Suppose $w\in S_{n}$ is compatible with
$\{a_1<a_2<\ldots< a_k\}$. Then we have
\begin{equation}
\label{eqn:splitA}
\schub_{w}(X)=\sum_{\underline \lambda}c_{\underline \lambda}(w)s_{\lambda^{1}}(X_1)\cdots
s_{\lambda^{k}}(X_{k})
\end{equation}
where $X_{i}=\{x_{a_{i-1}+1},\ldots,x_{a_{i}}\}$ and the sum is over all
sequences of partitions \linebreak ${\underline \lambda}=(\lambda^{1},\ldots,\lambda^{k})$. Each
$c_{\underline \lambda}(w)$ is a quiver coefficient, equal to the number of sequences
of semi-standard tableaux $(T_1,\ldots, T_{k})$ such that: 
\begin{itemize}
\item[(i)] $T_{1},T_{2},\ldots, T_{k}$ have entries strictly larger than $0,a_1,\ldots,a_{k-1}$ respectively;
\item[(ii)] the shape of
$T_i$ is conjugate to $\lambda^{i}$;
\item[(iii)] ${\rm col}(T_1)\cdots {\rm col}(T_k)$ is a
reduced word for $w$.
\end{itemize}
\end{thm}

        We will need some notation and definitions.
When $\mu=(\mu_1>\mu_2>\ldots>\mu_\ell)$
is a partition with $\ell$ distinct parts, there is a \emph{shifted shape}
given by a Ferrers shape of $\mu$ where each row is indented
one space from the left of the row above it. A \emph{shifted tableau} of 
shape $\mu$ is a filling of the shifted shape of $\mu$ by numbers and circled 
numbers $1^{\circ}<1<2^{\circ}<2<\ldots$ that is non-decreasing along each row and column. A shifted tableau is a {\em circled shifted tableau} if no 
circled number is repeated in any row and no uncircled number is repeated 
in any column. 
\begin{figure}[h]
\[\begin{picture}(120,100)
\put(0,70){\framebox(120,20)}
\put(20,70){\line(0,1){20}}
\put(40,70){\line(0,1){20}}
\put(60,70){\line(0,1){20}}
\put(80,70){\line(0,1){20}}
\put(100,70){\line(0,1){20}}
\put(20,50){\line(0,1){20}}
\put(40,50){\line(0,1){20}}
\put(60,50){\line(0,1){20}}
\put(80,50){\line(0,1){20}}
\put(20,50){\line(1,0){60}}
\put(40,30){\line(0,1){20}}
\put(60,30){\line(0,1){20}}
\put(80,30){\line(0,1){20}}
\put(40,30){\line(1,0){40}}
\put(60,10){\line(1,0){20}}
\put(60,10){\line(0,1){20}}
\put(80,10){\line(0,1){20}}
\put(6,76){$1^{\circ}$}
\put(26,76){$1$}
\put(46,76){$2$}
\put(66,76){$4$}
\put(86,76){$4$}
\put(106,76){$5^{\circ}$}

\put(26,56){$4^{\circ}$}
\put(46,56){$4$}
\put(66,56){$6^{\circ}$}

\put(46,36){$5^{\circ}$}
\put(66,36){$6^{\circ}$}

\put(66,16){$8$}

\end{picture}\]
\caption{A circled shifted tableau for $\mu=(6>3>2>1)$}
\end{figure}

The {\em weight} $x^{T}=x_1^{w_1} x_{2}^{w_2}\cdots$ of a 
circled shifted tableau is defined by setting $w_i$ to be the number of 
$i$ or $i^{\circ}$ occurring in $T$. With
this, the Schur $Q$ function $Q_{\mu}(X)$ is defined as $\sum_{T}x^{T}$, taken
over all circled shifted tableaux of shape $\mu$. The Schur $P$ function
$P_{\mu}(X)$ is defined to be $2^{-\ell(\mu)}Q_{\mu}(X)$, where $\ell(\mu)$ is the number
of parts of $\mu$ (see, e.g., \cite{fulton.pragacz,hoffman.humphreys}).      

        The Weyl group for the types $B_{n}$ and $C_{n}$ is the {\em hyperoctahedral group} ${\mathbb B}_{n}$ of
signed permutations on $\{1,2,\ldots,n\}$. It is generated by the simple
transpositions $s_{i}$ for $1\leq i\leq n-1$ together with the special
generator $s_0$, which changes the sign of the first entry of the
signed permutation.
The Weyl group of type $D_{n}$ is
the subgroup ${\mathbb D}_{n}$ of ${\mathbb B}_{n}$ whose elements make
an even number of sign changes. It is generated by the simple transpositions
$s_{i}$ for $1\leq i\leq n-1$ together with $s_{\hat 0}=s_0 s_1 s_0$. 
        
        The $B_n$ and $D_n$ analogues of Stanley functions,
$F_w(X)$ for $w\in {\mathbb B_{n}}$ and $E_{w}(X)$ for 
$w\in {\mathbb D_{n}}$, respectively, are defined in \cite{billey.haiman} by
\[F_{w}(X)=\sum_{\mu}f_{w\mu}Q_{\mu}(X)\]
and
\[E_{w}(X)=\sum_{\mu}e_{w\mu}P_{\mu}(X),\]
for certain nonnegative integers $f_{w\mu}$ and $e_{w\mu}$
given by explicit positive 
combinatorial formulas which we will not reproduce here; see \cite{billey.haiman} for details.

        In \cite{billey.haiman}, the theory of $A_{n-1}$ Schubert
          polynomials 
\cite{LS1} was extended to types $B_n,C_n$ and~$D_n$ (see
          \cite{fomin.kirillov:B_n} for an alternative approach). In each case, 
the corresponding generalized flag variety of order $n$ naturally projects into the one of order
$n+1$. This yields maps on the corresponding 
cohomology rings that sends Schubert classes to Schubert classes,
which in turn yields Schubert polynomials in the inverse limit. These are computed
as the unique solution of an infinite system of divided difference
equations. See \cite{billey.haiman} for details.

        For types $B_n,C_n$ and $D_n$,
the Schubert polynomials ${\mathfrak B}_{n}$, ${\mathfrak C}_{n}$ and
${\mathfrak D}_{n}$ respectively live in the polynomial ring
${\mathbb Q}[x_1,x_2,\ldots; p_{1}(Z),p_{2}(Z),\ldots]$, where $p_{k}(Z)=
z_{1}^k + z_{2}^k +\cdots$ is a power series in a new collection of variables
$Z=\{z_1,z_2,\ldots\}$. It is then proved in \cite{billey.haiman} that
for $w\in {\mathbb B}_{n}$,
\begin{equation}
\label{eqn:BH_C}
{\mathfrak C}_{w}=\sum_{u,v} F_{u}(Z)\schub_{v}(X),
\end{equation}
where the sum is over $u\in{\mathbb B}_{n}$ and $v\in S_{n}$ with $uv=w$ 
and $\ell(u)+\ell(v)=\ell(w)$. Also, if $s(w)$ is the number 
of sign changes of $w$, then
\begin{equation}
\label{eqn:BH_B}
{\mathfrak B}_{w}=2^{-s(w)}{\mathfrak C}_{w}.
\end{equation}
Similarly for $w\in {\mathbb D}_{n}$,
\begin{equation}
\label{eqn:BH_D} 
{\mathfrak D}_{w}=\sum_{u,v} E_{u}(Z)\schub_{v}(X),
\end{equation}
where the sum is over $u\in {\mathbb D}_n$ and $v\in S_{n}$, 
with $uv=w$ and $\ell(u)+\ell(v)=\ell(w)$. 

        More generally, if $w\in {\mathbb B}_{n}$ and a sequence of nonnegative integers $\{a_{j}\}$ with
$1\leq a_1<a_2<\ldots<a_k<n$, we say that $w$ is {\em compatible} with 
$\{a_{j}\}$ if whenever $\ell(ws_i)<\ell(w)$ for a simple transposition
$s_i$, then $i\in\{a_{j}\}$.

\begin{thm}
\label{thm:splitBCD}
Let $w\in {\mathbb B}_{n}$ be  
compatible with $\{a_1< a_2<\ldots< a_k\}$. Then we have
\begin{equation}
\label{eqn:splitC}
{\mathfrak C}_{w}=\sum_{\mu;{\underline \lambda}}c_{\mu;{\underline
    \lambda}}(w)
Q_{\mu}(Z)s_{\lambda^{1}}(X_1)s_{\lambda^2}(X_2)\cdots
s_{\lambda^{k}}(X_k)
\end{equation}
and
\begin{equation}
\label{eqn:splitB}
{\mathfrak B}_{w}=2^{-s(w)}\sum_{\mu;{\underline
    \lambda}}c_{\mu;{\underline \lambda}}(w)
Q_{\mu}(Z)s_{\lambda^{1}}(X_1)s_{\lambda^2}(X_2)\cdots
s_{\lambda^{k}}(X_k).
\end{equation}
If in addition, $w\in {\mathbb D}_{n}$, then
\begin{equation}
\label{eqn:splitD}
{\mathfrak D}_{w}=\sum_{\mu;{\underline \lambda}}d_{\mu;{\underline
    \lambda}}(w)
P_{\mu}(Z)s_{\lambda^{1}}(X_1)s_{\lambda^2}(X_2)\cdots
s_{\lambda^{k}}(X_k).
\end{equation}
In the above formulas, $X_{i}=\{x_{a_{i-1}+1},\ldots, x_{a_{i}}\}$,
$\mu$ is a partition with distinct parts and ${\underline
  \lambda}=(\lambda^1,\ldots,\lambda^k)$ is a sequence of
partitions. Also, 
$c_{\mu;{\underline\lambda}}(w)=f_{u\mu}c_{\underline \lambda}(v)$ and
$d_{\mu;{\underline \lambda}}=e_{u\mu}c_{\underline \lambda}(v)$ where
$uv=w$, $\ell(u)+\ell(v)=\ell(w)$, $v\in S_{n}$, and $u\in {\mathbb
  B_{n}}$ or $u\in {\mathbb D}_{n}$, respectively.
\end{thm}
\begin{proof}
Suppose $w\in {\mathbb B_{n}}$ (or respectively, $w\in {\mathbb
  D}_{n}$) and $uv=w$ with $\ell(u)+\ell(v)=\ell(w)$ where $u\in
  {\mathbb B_{n}}$ (or $u\in {\mathbb
  D}_{n}$) and $v\in S_{n}$.

Let $i\geq 1$ be such that $\ell(vs_i)<\ell(v)$. Then by our assumptions
and
standard properties of the length function (see, e.g., \cite[Section~5.2]{Humph}) we have
\[\ell(ws_i)=\ell(uvs_i)\leq \ell(u)+\ell(vs_i) < \ell(u)+\ell(v)=\ell(w).\]
Hence $i$ is one of the $a_{j}$, i.e., $v$ is compatible with $\{a_{j}\}$.  
Therefore, the result follows from equations (\ref{eqn:BH_C}),(\ref{eqn:BH_B}) and $(\ref{eqn:BH_D})$ combined with Theorem~\ref{thm:Schub_split}.
\end{proof}

\begin{exa}{\rm  Consider $w=\left(
\begin{array}{ccccc}
1 & 2  & \ \ 3 \\
3 & 1  & -2 
\end{array}
\right)=s_1 s_0 s_1 s_2 s_1\in {\mathbb B}_{3}$. 
This signed permutation is compatible with the
    sequence $1<2$. In \cite{billey.haiman} the following was
    computed:
\[{\mathfrak C}_{w}=Q_{41}+Q_4 x_1 + Q_{31}x_1 + Q_3 x_1^2 + Q_{31}x_2
    + Q_{3}x_1 x_2 + Q_{21}x_1 x_2 + Q_2 x_1^2 x_2.\]
This may be rewritten as
\begin{multline}
{\mathfrak C}_{w}=Q_{41} + Q_{4}s_{1}(x_1) +Q_{31}s_{1}(x_1)
+Q_{31}s_{1}(x_2) + Q_{3}s_{2}(x_1) + \\
Q_{3}s_{1}(x_1)s_{1}(x_2)
+Q_{21}s_{1}(x_1)s_{1}(x_2) +Q_{2}s_{2}(x_1)s_{1}(x_2),
\end{multline}
in agreement with Theorem~\ref{thm:splitBCD}.
}
\end{exa}
        In \cite{BKTY} it was explained why (\ref{eqn:splitA})
provides a geometrically natural solution to the {\em Giambelli} problem
for partial flag varieties. For the other classical types,
the choice of variables makes it unclear what the underlying
geometry of (\ref{eqn:splitC}),(\ref{eqn:splitB}) and
$(\ref{eqn:splitD})$ might be. On the other hand, given the shape of the
formulas, by analogy with the $A_{n-1}$ case,  it is natural to ask if
there is a degeneracy locus setting for which the coefficients 
$c_{\mu;{\underline \lambda}}(w)$ and $d_{\mu;{\underline
    \lambda}}(w)$ (and their positivity)
appear.

\section{Proof of Theorem~\ref{thm:main_1}}
Let $S_{d}(\rr)$ denote the set of permutations $w$ in $S_{d}$ such that
$G(w)$ contains the same number of points in $M_{ij}$ as $G(v(\rr))$
does, for all $0\leq i,j\leq n$. Our proof of Theorem~\ref{thm:main_1}
uses the following:
\begin{prop}
\label{thm:characterization}
Let $\rr=\{r_{ij}\}$ for $0\leq i\leq j\leq n$ be a set of rank
conditions and let $w\in S_{d}$, $d=r_{00}+r_{11}+\cdots +
r_{nn}$. The following are equivalent:
\begin{itemize}
\item[(I)] $w=v(\rr)$;
\item[(II)] $w$ is the minimal length element of $S_{d}(\rr)$;
\item[(III)] there exists a pipe dream $D$ for $w$
  and there exists a lacing diagram $W$ for $\rr$ such that 
$D$ has every box in $\bigcup_{i+j\leq n-2} M_{ij}$ tiled by crosses, 
$D$ maps to $W$, and $\ell(w)\leq \ell(v(\rr))$.  
\end{itemize}
\end{prop}
\begin{proof}
The length of $w\in S_{d}(\rr)$ is computed from $G(w)$ by counting
those pairs of dots where one is situated to the northeast of the
other. Call such a pair {\em unavoidable} if the dots actually appear
in blocks where one is situated (strictly) northeast of the
other. The number of unavoidable pairs is constant on
$S_{d}(\rr)$. Moreover, observe that all of the pairs contributing to the length of
$v(\rr)$ are unavoidable. On the other hand, if $w\neq v(\rr)$, then at
least one pair contributing to $\ell(w)$ is not unavoidable. Thus 
(I) is equivalent to (II).

That (I) implies (III) is immediate from
\cite[Theorem~5.10]{KMS}, but we include a proof for
completeness. Take any $D\in {\mathcal RC(v(\rr))}$. The definition of
$v(\rr)$ implies that $D$ has every box in $\bigcup_{i+j\leq n-2} M_{ij}$ tiled by crosses.
Moreover, $D$ gives a
lacing diagram $W$ such that $D$ maps to $W$. Observe that
the number of pipes of $D$ that enter in the $i^{th}$ horizontal strip
and exit in the $j^{th}$ vertical strip is equal to the number of
points of $G(v(\rr))$ in $M_{ij}$ for $0\leq i,j\leq n$. From this and
the definition of $v(\rr)$ it
follows that $W$ is in fact a lacing diagram for $\rr$. 

Finally, suppose (III) holds. 
By considering where the pipes of $D$ go in relation to $W$, one finds
that $G(w)$ and $G(v(\rr))$ have the same number
of points in any block on the main anti-diagonal and below, i.e.,
blocks
$M_{ij}$ where $i+j\geq n$. The 
condition on the boxes of $\bigcup_{i+j\leq n-2}M_{ij}$ implies that the
only other points of $G(w)$ appear in
the blocks $M_{i,n-i-1}$, $0\leq i\leq n-1$ on the super-antidiagonal. Since $w$ is a permutation,
each of these blocks must have the same number of points as its counterpart in
$G(v(\rr))$, i.e. $w\in S_{d}(\rr)$. Since we already know $v(\rr)$ is
the unique minimal length element of $S_{d}(\rr)$, the assumption that
$\ell(w)\leq \ell(v(\rr))$ implies~(II).
\end{proof}

\subsection{Proof of Theorem~\ref{thm:main_1}}
Let $\rho$ be a partial permutation represented by an $a\times b$ matrix. Consider the diagram $D({\tilde \rho})$
of ${\tilde \rho}$, which consists of the boxes $(i,j)$ in $(a+b)\times (a+b)$ such that
${\tilde \rho}(i)>j$ and ${\tilde \rho}^{-1}(j)>i$. 
Associated to ${\tilde \rho}$ is its {\em canonical reduced word}. This
 is obtained by numbering the boxes of $D({\tilde \rho})$ consecutively in each row,
from right to left, starting with the number of the row. Then the rows are read
left to right, from top to bottom (see, e.g., \cite{Manivel}).

\begin{lemma}
\label{lemma:exists_red}
Let $u_1 u_2 \cdots u_{\ell({\tilde \rho})}$ be the
canonical reduced word for ${\tilde \rho}$. 
Then the set $\{k_1<k_2<\ldots<k_p\}$ of indices $k$ where
$u_k<u_{k+1}$
has size at most $a$. Moreover,
$j\leq u_{k}$ for all $k\in [k_{j-1}+1,k_j]$, where $k_{0}=0$.
\end{lemma}
\begin{proof}
By construction, $D({\tilde \rho})$ sits inside the northwest $a\times
b$ rectangle of the \linebreak $(a+b)\times (a+b)$ box. Since the labels of
the boxes in the construction of the canonical reduced word decrease
from left to right along each row, there can be at most $a$ indices
$k$ where $u_{k} < u_{k+1}$. The fact that each entry of the $t^{th}$
row of the filling of $D({\tilde \rho})$ is at least $t$ implies the
remainder of the claim.
\end{proof}

\begin{exa}{\rm
Let $\rho$ be the partial permutation represented by the matrix:
\[\left(\begin{matrix}
0 & 0 & 0 & 0 \\
1 & 0 & 0 & 0 \\
0 & 0 & 0 & 1 
\end{matrix}\right).\]
The canonical reduced word for ${\tilde \rho}$ is obtained below (see Figure~2).

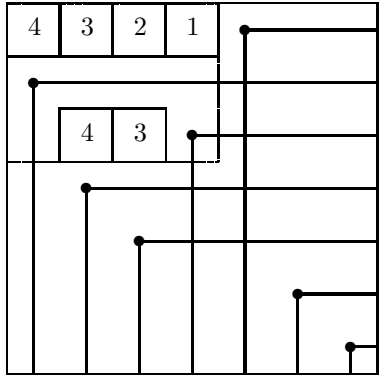
\begin{figure}[h]
\label{fig:canonical_red_word}
\[
\begin{picture}(140,140)
\put(0,0){\framebox(140,140)}
\put(90,130){\circle*{4}}
\thicklines
\put(90,130){\line(1,0){50}}
\put(90,130){\line(0,-1){130}}
\put(10,110){\circle*{4}}
\put(10,110){\line(1,0){130}}
\put(10,110){\line(0,-1){110}}

\put(70,90){\circle*{4}}
\put(70,90){\line(1,0){70}}
\put(70,90){\line(0,-1){90}}

\put(30,70){\circle*{4}}
\put(30,70){\line(1,0){110}}
\put(30,70){\line(0,-1){70}}

\put(50,50){\circle*{4}}
\put(50,50){\line(1,0){90}}
\put(50,50){\line(0,-1){50}}

\put(110,30){\circle*{4}}
\put(110,30){\line(1,0){30}}
\put(110,30){\line(0,-1){30}}

\put(130,10){\circle*{4}}
\put(130,10){\line(1,0){10}}
\put(130,10){\line(0,-1){10}}

\thinlines
\put(0,80){\dashbox{.5}(80,60)}

\put(0,120){\framebox(20,20)}
\put(20,120){\framebox(20,20)}
\put(40,120){\framebox(20,20)}
\put(60,120){\framebox(20,20)}

\put(8,128){$4$}
\put(28,128){$3$}
\put(48,128){$2$}
\put(68,128){$1$}

\put(20,80){\framebox(20,20)}
\put(40,80){\framebox(20,20)}

\put(28,88){$4$}
\put(48,88){$3$}
\end{picture}
\]
\caption{The canonical reduced word $4321\cdot 43$ for ${\tilde \rho}$}
\end{figure}
}
\end{exa}

The following fact is immediate from the main theorem of \cite{BB}. We include
a proof for completeness:
\begin{lemma}
\label{prop:local}
There exists an RC-graph for ${\tilde \rho}$ such that
all crosses occur in its northwest $a\times b$ sub-rectangle.
\end{lemma} 
\begin{proof}
Let $u_{1} u_2 \cdots u_{\ell({\tilde \rho})}$ be the canonical reduced
word for ${\tilde \rho}$. By Lemma~\ref{lemma:exists_red},
\begin{equation}
\label{eqn:rcs}
(\underbrace{1,1,\ldots,1}_{k_{1}},\underbrace{2,2,\ldots,2}_{k_{2}},
\ldots,\underbrace{p,p,\ldots,p}_{k_{p}})
\end{equation}
is a reduced compatible sequence for ${\tilde
  \rho}$, and the conclusion follows from Proposition~\ref{prop:rctoRC}.
\end{proof}

\begin{exa}{\rm Continuing the previous example, the reduced
    compatible sequence $(\ref{eqn:rcs})$ corresponding to the
canonical reduced word for ${\tilde \rho}$ is 
\[(1,1,1,1,2,2).\]
By Proposition~\ref{prop:rctoRC}, there is an RC-graph for ${\tilde
  \rho}$ with crosses from
\[\{(1,4),(1,3),(1,2),(1,1),(2,3),(2,2)\}.\]
That RC-graph is 
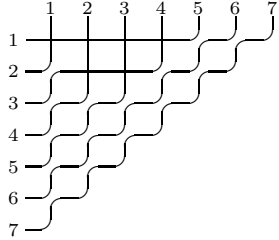
\begin{figure}[h]
\[
\begin{array}{ccccccccc}
       &\perm1{}&\perm2{}&\perm3{}&\perm4{} & \perm5{} & \perm6{} & \perm7{} \\
\petit1 & \+ & \+ & \+  & \+ & \jr & \jr  & \je  \\
\petit2 & \jr & \+ & \+ & \jr & \jr & \je  &     \\
\petit3 & \jr & \jr & \jr & \jr & \je &    &     &      \\
\petit4 & \jr & \jr & \jr & \je &     &     &     &      \\
\petit5 & \jr & \jr & \je &  & &     &     &      \\
\petit6 & \jr & \je & &  & &     &     &      \\
\petit7 & \je & &  &  & &     &     &      \\
\end{array}
\]
\end{figure}
}
\end{exa}



\noindent
\emph{Proof of Theorem~\ref{thm:main_1}.}
Construct a pipe dream $D$ starting with a $d\times d$ box as follows. 
For $k=1,2,\ldots,n$ let $D_{k}$ be the RC-graph obtained by applying 
Lemma~\ref{prop:local} to the partial permutation $w_{k}$. Then let
${\overline D_{k}}$ denote the northwest $r_{k-1,k-1}\times r_{kk}$
sub-pipe dream, rotated 180 degrees. Overlay
${\overline D_{k}}$ into $M_{k-1,n-k}$. 
For the remaining boxes, place crosses
in the top $r_{00}+r_{11}+\cdots +r_{n-2,n-2}$
rows of the $d\times d$ box and elbow joints elsewhere. This defines
a pipe dream $D$ for some permutation $w\in S_{d}$. 
\begin{figure}[h]
\[
\begin{picture}(150,150)
\put(0,0){\framebox(150,150)}
\put(0,0){\dashbox{.5}(30,150)}
\put(0,0){\dashbox{.5}(70,150)}
\put(0,0){\dashbox{.5}(100,150)}
\put(0,0){\dashbox{.5}(120,150)}
\put(0,0){\dashbox{.5}(150,30)}
\put(0,0){\dashbox{.5}(150,70)}
\put(0,0){\dashbox{.5}(150,100)}
\put(0,0){\dashbox{.5}(150,120)}
\thicklines
\put(0,30){\framebox(30,40)}
\put(30,70){\framebox(40,30)}
\put(70,100){\framebox(30,20)}
\put(100,120){\framebox(20,30)}
\put(11,48){${\overline D_{n}}$}
\put(45,85){$\cdots$}
\put(78,107){${\overline D_{2}}$}
\put(105,130){${\overline D_{1}}$}
\thinlines
\put(0,145){\line(1,0){100}}
\put(0,135){\line(1,0){100}}
\put(0,125){\line(1,0){100}}
\put(0,115){\line(1,0){70}}
\put(0,105){\line(1,0){70}}
\put(0,95){\line(1,0){30}}
\put(0,85){\line(1,0){30}}
\put(0,75){\line(1,0){30}}

\put(5,150){\line(0,-1){80}}
\put(15,150){\line(0,-1){80}}
\put(25,150){\line(0,-1){80}}
\put(35,150){\line(0,-1){50}}
\put(45,150){\line(0,-1){50}}
\put(55,150){\line(0,-1){50}}
\put(65,150){\line(0,-1){50}}
\put(75,150){\line(0,-1){30}}
\put(85,150){\line(0,-1){30}}
\put(95,150){\line(0,-1){30}}
\end{picture}\]
\caption{Construction of $D$}
\end{figure}
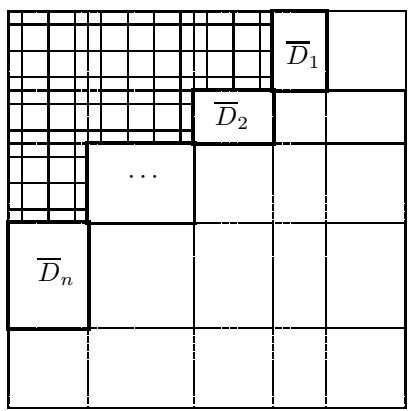
By construction, $D$ maps to $W$ and 
moreover, the number of crosses in $D$ is
\[\mid\!\bigcup_{i+j\leq n-2} M_{ij}\!\mid + \ell(W).\]
Since $W$ is minimal length, $\ell(W)=d(\rr)$ and so
by (\ref{eqn:crucial}), $l(w)\leq l(v(\rr))$. Then by
Proposition~\ref{thm:characterization}, $w=v(\rr)$ and thus $D\in
{\mathcal RC}(v(\rr))$.
This construction describes the desired injection.\qed
\medskip

        For example, the RC-graph given in Figure~1 is the image of $W$ from
Example~1 under the embedding map of Theorem~\ref{thm:main_1}.




\end{document}